# Inversion of finite filters

**YURY TYURIN AND ANASTASIA TYURINA**
**SECOND STAR ALGONUMERIX LLC**
**MAY 2022**

**Abstract**
We present a decomposition of finitely supported filters ( aka instrument function PSF) as a composition of invertible and non-invertible filters. The invertible component can be inverted directly and the non-invertible component is shown to decrease the resolution of the acquired signal.
**Summary**
The effects of instruments of data acquisition on reality are typically described by 'running window' operators also called convolutions and finite filters. We consider operations of convolutions with finite symmetrical filters acting on finite sequences. We introduce characteristic polynomials for symmetric filters. Characteristic polynomials of convolutions are products of the characteristic polynomials of the convolved filters. We show that as polynomials of even degree allow for decomposition into quadratic polynomials with real coefficients, finite symmetrical convolutions allow for a decomposition into elementary symmetric filters of length three. We show the conditions for the elementary filters to be invertible. Any filter is a composition of invertible and non-invertible elementary filters and can be presented as a convolution of an invertible and a non-invertible components . The invertible component can be undone via convolution with its inverse. Each non-invertible elementary convolution of the non-invertible component leads to loss of resolution in the signal (compared to the original sequence). The problem allows for a 2D generalization for separable filters.

# 1. Statement of the problem

Our instruments measure certain characteristics of objects (height, incline, color etc) and report the observed values of object under measurement. In this paper we will start with one dimensional signals $X = \{ x(t), t \in Z)\}$ parametrized by coordinate t. The value $Y = \{ y(t), t \in Z )\}$, reported by the instrument is a linear combination of the values in the neighboring points:

$$y(t) = \sum c(k) x(t-k) \qquad (1.1)$$

Here C is a finite and (typically) symmetrical sequence on integers $C = \{c(t), (t \in Z)\}$ referred to as a filter, and it is a fundamental characteristic of a measurement instrument. It is called Point Spread Function (PSF). The statement (1.1) can be written as a convolution of X and the filter C:

$$Y = C*X \qquad (1.2)$$

In this work the PSF of the instrument is considered "known".

The authors are expressing their gratitude to Dr. Yaschouk for his discussions that got us interested in the general problem of inversion of finite symmetrical filters.

# 2. Inverse convolution

Let $I = \{I(t) / t \in Z\}$ be a *unitary* sequence, such that $I(0) = 1$, $I(t) = 0$ for any t, not equal to 0. Then a unitary sequence is a unitary filter, so that a convolution with it change no sequence X:

$$X = I*X \qquad (2.1)$$

If *C* is a given finite sequence of order N,

$$C: k \to c(k), \text{ for } |k| \leqslant N \text{ and } C: k \to 0, \text{ for } |k| \geqslant N \qquad (2.2)$$

Add the condition of symmetry of the filter *C*: $c(-k) = c(k)$, for any $k$, $|k| \leqslant N$ $\qquad (2.4)$

Let sequence $Z = \{z(t) \in Z\}$, be an inverse to *C* if $\quad C * Z = I,$ $\qquad (2.5)$

If Z is found we can solve (1.2) for *X* as

$$X = Z * Y \qquad (2.6)$$

We should state that not every sequence has an inverse. A serviceable inverse sequence, should be summable and its elements should converge to zero exponentially fast (when $t \to \infty$). If a filter-sequence *C* has an inverse we will call it ***invertible.***

# 3. Characteristic polynomial of a finite convolution

Let us introduce characteristic polynomials of convolutions. A finite filter *C* introduced in (1.1) and (2.2) is acting by a convolution with a finite sequence *C* of order *N* as defined in (2.3). We define its characteristic polynomial *P* (of order *2N*) as :

$$P_{2N}(x) = \sum_{k=-N}^{k=N} c(k) x^k \qquad (3.1)$$

If the filter *C* is symmetric and $c(-N) = c(N) \neq 0$, the polynomial can be normalized by $c(N)/x^N$

$$P_{2N}(x) = x^{2N} + \sum_{k=1}^{k=2*N-1} c(k)/c(N) * x^k + 1 \qquad (3.2)$$

A characteristic polynomial of a symmetric sequence of order *N* is a symmetric polynomial of order *2N* with real coefficients. Consider a factorization of *P*:

$$P_{2N}(x) = \prod (x - x(l)), \qquad (3.3)$$

here $x(l)$ (l = 1, …, 2N) are all the roots of the polynomial $P_N$.

As $P_{2N}$ is symmetric, its roots come in pairs. Namely, for a root *u*, its inverse 1/u is also a root of $P_{2N}$. In each such pair of the real roots $(u, 1/u)$ the absolute value one of the roots is greater than one and that of the other is less than one. A dual pair of roots may also consist of two imaginary numbers u and its conjugate, which is in the same time $\frac{1}{u}$. Therefore, if the polynomial P has any complex roots, they are all of the magnitude one. Then the polynomial $P_N(x)$ can be presented as

$$P_N(u) = \prod (u^2 - p(k)u + 1), \text{ where } p(1), \ldots, p(N) - \text{are real numbers.} \qquad (3.4)$$

For a convolution of two sequences its characteristic polynomial is a product of their respective characteristic polynomials:

$$P_{l+k}(C_l * C_k) = P_l(C_l) P_k(C_k) \qquad (3.5)$$

If we consider a filter-sequence $C_m$ and its characteristic polynomial $P(C_m)$ and present the polynomial as a product of lesser degree polynomials as in Equation (2.1.5) we reduced the degree of the sequences that we are seeking to inverse. Namely, if we manage to inverse $C_k$ and $C_l$ separately then the inverse to their convolution is the convolution of the inverses. Let's reduce the degree as far as we can.

## 4. Elementary convolutions aka elementary symmetric filters

A quadratic polynomial $(u^2 + p(k)u + 1)$ is a characteristic polynomial for the symmetric filter sequence *(1, p(k), 1)*–of the length 3 and order 1:

$$c(k) = 0, for |k| > 1, c(-1) = 1, c(0) = p(k), c(1) = 1 \qquad (4.1)$$

We will call symmetric sequences of length 3 (and order 1) ***elementary symmetric filters***. We call them of order one as each them is defined by just one parameter *p(k)*.

Let's call the filter sequence C (introduced in 1.2), $C_N$ in order to specify the order of *C* (introduced in 1.2.3). The elementary sequences are of order 1 and should be denoted as $C_1$. Every $C_N$ is a convolution of *N* elementary sequences:

$$C_N = C_1^1 * C_1^2 * \ldots * C_1^N \qquad (4.2)$$

The *sequence* $C_N$ *(4.2) is invertible if and only if each of the elementary sequences* $C_1^1, C_1^2, \ldots, C_1^N$ *is invertible. If the inverse of* $C_N$ *exists, we can construct it explicitly.*

# 5. Decomposition of convolutions to compositions of elementary convolutions

How can we find the decomposition of a symmetrical filter $C_N$ into a convolution of elementary filters? We need to find a decomposition of $P_{2N}$, the characteristic polynomial of $C$ into a product of $N$ quadratic polynomials:

$P_{2N} = P_1^1 P_1^2 \ldots P_1^N$  (5.1)

There a simple method to find the quadratic decomposition for symmetric polynomials (3.2). To find the coefficients $p(1), \ldots, p(k), \ldots, p(N)$ we introduce a new polynomial $Q_N(x) = \prod_{k=1}^{k=N}(x - p(k))$ with roots equal to the coefficients.

The coefficients of the polynomial are the elementary symmetric polynomials $\sigma(p(1), \ldots, p(k), \ldots, p(N))$.

$$Q_N(x) = \sum_{k=0}^{k=N}(-1)^k \sigma_{(N-k)}((p(1),\ldots,p(N)))x^k \quad (5.2)$$

Here $\sigma_n(p(1), \ldots, p(N))$ is a sum of monomials of order $n$ of $(p(1), \ldots, p(N))$, namely
$\sigma_0(p(1), \ldots, p(N)) = 1$
$\sigma_1(p(1), \ldots, p(N)) = p(1) + p(2) + \ldots + p(N)$,
$\sigma_2(p(1), \ldots, p(N)) = p(1)p(2) + \ldots + p(1)p(N) + \ldots + p(N-1)p(N) \ldots$
$\sigma_N(p(1), \ldots, p(N)) = p(1)\ldots p(N)$

$\sigma_1, \ldots, \sigma_m$ can be expressed in terms of the coefficients of the characteristic polynomial $P_N$ (3.2).

**Example N=2**.
$P_2(x) = x^4 + c(1)x^3 + c(0)x^2 + c(1)x + 1 = (x^2 + p(1)x + 1)(x^2 + p(2)x + 1)$
Therefore $\sigma_1(p(1), p(2)) = c(1)$, $\sigma_2(p(1), p(2)) = c(0) - 2$ and
$Q_2(x) = x^2 - c(1)x + c(0) - 2$.

**Example N= 3**.
$P_3(x) = x^6 + c(2)x^5 + c(1)x^4 + c(0)x^3 + c(1)x^2 + c(2)x + 1$,
$P_3(x) = (x^2 + p(1)x + 1)(x^2 + p(2)x + 1)(x^2 + p(3)x + 1) = x^6 + [p(1) + p(2) + p(3)]x^5 +$
$+ [p(1)p(2) + p(1)p(3) + p(2)p(3) + 3]x^4 + [2(p(1) + p(2) + p(3)) + p(1)p(2)p(3)]x^3 + \ldots$
Equating the coefficients for the same powers of x:
$\sigma_1(p(1), p(2), p(3)) = c(2)$,
$\sigma_2(p(1), p(2), p(3)) = c(1) - 3$,
$\sigma_3((p(1), p(2), p(3)) = c(0) - 2c(2)$.
Therefore $Q_3(x) = x^3 - c(2)x^2 + (c(1) - 3)x - (c(0) - 2c(2))$.

For $N$ larger than 3 the method proceeds similarly and results in a linear system of equations connecting polynomials $\sigma_1 \ldots, \sigma_N$ and the known coefficients of the filter $C$ of any order $N$. The iterative procedure is closely related to the iterative process called Newton's identities using Chebyshev polynomials.

# 6. Inversion of elementary convolutions

Let us construct an inverse sequence $Z, Z = \{z(t) \in Z\}$ of an elementary convolution $C_0$.

We have to solve $C_0 * Z = I$                                                                           (6.1)
where *I* is a unitary convolution defined in (3.1)
For the elementary filer $C_0 = [1, p, 1]$ we have $z(t-1) + p * z(t) + z(t+1) = I(t)$     (6.2)

If we go through the consecutive values of $t: t=0, t=1, \ldots ut=-1, \ldots$, we will obtain linear equations connecting the values of elements of *Z(t)*:
1. $t=0$:   $z(-1) + p * z(0) + z(1) = 1 \quad z(1) = 1 - p * z(0) - z(-1)$
2. $t=1$:   $z(0) + p * z(1) + z(2) = 0 \quad z(2) = -p * z(1) - z(0)$
3. $t=2$:   $z(1) + p * z(2) + z(3) = 0 \Rightarrow z(3) = -p * z(2) - z(1)$

and so on, for any $t: z(t+1) = -p * z(t) - z(t-1)$                                    (6.3)
If we assume Z to be symmetrical in a sense that $z(-k) = z(k)$ we can only expend in positive direction and in the equation 1) we can assume $z(-1) = z(1)$ and solve it for $z(1)$:
$z(1) = ½ * (1 - p * z(0))$                                                              (6.4)
Let us define a vector $\vartheta(n) = (z(n), z(n-1))$, then $\vartheta(n+1) = \vartheta(n) * B$, where B is a linear operator and represented by a matrix $B = \begin{bmatrix} -p & 1 \\ -1 & 0 \end{bmatrix}$.                                             (6.5)

It follows that
and $\vartheta(1) = (z(1), z(0))$                                                (6.6)

**Asymptotic properties of $\vartheta(n)$ with the growth of n**

As we stated we are only interested in the summable sequences *Z(t)*. Better yet if the elements of *Z* converge to zero exponentially fast. To understand the asymptotic behavior of *Z* we study understand that of $\vartheta(n)$. In fact, the asymptotic behavior of $\vartheta(n)$ depends on the eigenvalues of matrix *B* (described in 3.1/5). The eigenvalues are the roots of its characteristic polynomial: $det(B - u * I)$.
It is a familiar polynomial in u.
$det(B - u * I) = u^2 + pu + 1$                                                    (6.7)
By design the polynomial is the characteristic polynomial of the elementary convolution $C_0$.
Let's see how the asymptotic behavior depends on p:
**Case 1**. If $p<2$, then the polynomial has two distinct real roots: $u1$ and $u2$. One of them is greater than one and the other lees than 1 $(u1<1; u2>1)$ because their product is 1. The eigenvector of B, corresponding to $u1$ is proportional to a vector $(u1,1)$. In that case we can use as a starting point
$\vartheta(1) = (u1 * z(0), z(0))$                                                (6.8)
We have to choose z(0) such that the triplet satisfies the initial condition given by the equation (1): $2z(1) + z(0) = 1$ and that is possible.
If $z(0)$ is chosen this way we observe that
$\vartheta(n+1) = u_1 * \vartheta(n) = u_1^n * \vartheta(1)$                                       (6.9)
Converges to zero exponentially fast, as u1 was chosen to be less than 1.

***Therefore, if and only if $p>2$ the elementary convolution $C_1$, $[1, p, 1]$ is invertible and the inverse can be easily constructed explicitly. In fact, Z(t) converges to zero so fast that for all computational effects and purposes only a few elements of the filter Z have non-zero values (the filter has a finite support).***

**Case 2**. If $|p|=2$ the eigenvalues of B are multiples and equal to either $+1$ or $-1$. In that case $\vartheta(n)$ will not converge to 0 regardless of the choice of $\vartheta(1)$.

**Case 3.** If $|p|<2$ the eigenvalues of B are 'unitary' $|u_1|=1; |u_2|=1$,- they are complex, conjugate of each other. and their absolute value equal to 1. In fact, the two roots $u_1$ and $u_2$ are connected by Vieta equations:

$$u_1 * u_2 = 1; u_1 + u_2 = -p \tag{6.10}$$

Because $\|u1\|=1$, $\vartheta(n)$ will not converge to 0 regardless of the choice of $\vartheta(1)$.

***Therefore if $|p| \geq 2$ the elementary convolution $C_0$ ,$([1,p,1])$ is not invertible in terms of convergence to 0 as t tends to ∞. However, a solution that does not converge to zero (as t tends to ∞) does exist and can be constructed in close terms. We call this solution 'pseudoinverse'. The example of such a pseudoinverse for p = 1 is shown in Drawing 2c, and for p= - 1 in the Drawing 2d. The pseudo inverses oscillate: the elements of the filter do not converge to zero, but their values are bounded and oscillate between the numbers -0.5 and 0.5.***

# 7. Kernels of non-invertible elementary filters

Let's concentrate on non-invertible convolutions.

If $C_0$ is not invertible there should be a kernel $K(C_0)$ to the operator

$F \to C_0 * F$ , such that if $f \in K : f * C_0 = 0$ (7.1)

We offer ***an explicit method of construction of the kernel on a non-invertible elementary convolution***.

The kernel in our model is a discreet sequence of the length N (the length of the original data segment X) and each coordinate of it is given by the formula equation:

$X(n-1) + pX(n) + X(n+1) = 0$ (7.2)

The iterative procedure is somewhat similar to the method of construction of the inverse. Let's start with n = 0:

$X(-1) + p * X(0) + X(+1) = 0$ (7.3)

Assume symmetric Kernel and $X(-1) = X(+1)$, then $p * X(0) + 2 * X(1) = 0$ and $X(1) = -p/2 * X(0)$

For n = 1: $X(0) + p * X(1) + X(2) = 0$ and $X(2) = -p/2 * X(0)$ (7.4)

Denote a vector $\vartheta(n) = \{X(n), X(n-1)\}$, and $\vartheta(n+1) = B * \vartheta(n)$, then matrix $B = \begin{bmatrix} -p & -1 \\ 1 & -0 \end{bmatrix}$. Here $|p|=2$ as C is non-invertible.

It follows that

$\vartheta(n) = \vartheta(n-1) * B$ and $\vartheta(1) = (X(1), X(0))$ (7.5)

To find the eigenvalues of the matrix (and its eigenvectors) we consider its characteristic polynomial

$det\begin{pmatrix} -p-\lambda & 0 \\ 1 & -\lambda \end{pmatrix} = \lambda^2 + p\lambda + 1$ and the eigenvalues of the matrix are $\frac{-p}{2} \mp \frac{\sqrt{(p^2-4)}}{2}$ (7.6)

The recursive formula for the elements of the kernel $X(n)$ easily follows and allows us to evaluate it in close terms.

***In a decomposition of a PSF filter into elementary convolutions we identify the non-invertible elementary convolutions and their composition forms the non-invertible part of the PSF filter. For this non-invertible part of PSF we identify the kernel and form a factor space over it. The***

*kernel of each non-invertible elementary convolution has a dimension of 2. Therefore each non-invertible elementary convolution in the decomposition of the PSF reduces the length of the finite signal by 2.*

**Non-invertible filters and loss of resolution**

Nyquist theorem states that any (periodic and continuous on [-π, π]) signal *X* that does not have frequencies higher than that of cos(*Nx*) *and sign(Nx)* in its Fourier representation can be presented without loss of accuracy by its values in *2N+1* points (on [-1, 1]). That is because if *X* allows for Fourier transform *F(X)* (and continuous signals do), then

$$F(X) = \sum_{n=-\infty}^{n=+\infty} a_n \cos(nt) + b_n \sin(nt) \qquad (7.8)$$

and if $a_n = 0; b_n = 0, if\ n > N$ the equation (7.8) becomes

$$F(X) = \sum_{n=-N}^{n=+N} a_n \cos(nt) + b_n \sin(nt) \qquad (7.9)$$

In the Equation (4.3.1) there are 2N+1 unknowns $\{a_n, b_n\}, 0 < n < N$ and to find them we need to form *2N+1* equations, and for that knowing the values of *X* in exactly *2N* points would be required.

Suppose the interval, between the consecutive points in the signal representation was selected to be exactly the best estimate for the Nyquist frequency of the system. (In all likelihood the known Nyquist frequency was eliminated without loss of generality in the preprocessing of the signal and in fact 't' is simply measured in the units of length equal to Nyquist frequency and thus in the previous chapters we considered the signal represented on the lattice of [*1, 2, …, t, N*]).

If the acquired signal Y, (*Y = X\*C*) does not have the highest frequencies of cos(*Nx*) and sin(*Nx*) in its representation, because their convolution with *C0* is 0. That means that Y only have frequencies not higher than *(N-1)* in its representation and that means that the Nyquist interval required for its lossless sampling is $\frac{1}{2(N-1)} = \frac{1}{2N-2}$. That means that convolution with *C0* changed the Nyquist interval (or resolution) from $\frac{1}{2N}$ to $\frac{1}{2N-2}$. It is consistent with the notion that a non-invertible elementary convolution has a kernel space of dimension 2.

One may return to the "primitive argument" of solving (1.2) via Fourier transform. Fourier transform both sides of the equation, and the operation of convolution becomes multiplication. *Fourier(Y) = Fourier( C)\*Fourier (X) => Fourier(X) = Fourier(Y)/Fourier(C )) iff  Fourier(C ) # 0.*

The approach meets an obstacle if the Fourier transform of C has zero values.

When we identify the kernel of the convolution that corresponds to the zeros of the PSF in the Fourier domain. When we combine factorization over the kernel the operation reduces the size of the image, but very slightly - one point for every kernel function.

**Why are we able to invert some elementary convolutions then?**

Because, the Fourier transforms of elementary convolutions with |p|>2 do not have zero values.

## 8. Generalization to separable 2D filters

We can generalize the 1D method we presented for a 2D:

$Y(s,t) = C(s,t) *X(s, t)$                                                  (8.1)

We assume that the two-dimensional instrument function *C* (MTF) is symmetrical and compactly supported in each direction (limited to a finite rectangle). Consider a typical example where the instrument function is an inverse of the exponent of the square distance to the center of the running window within the window W:

$C = \exp(-s^2 - t^2) = \exp(-s^2)\exp(-t^2)$ for $(s,t) \in W$         (8.2)

In that case the 2D function C is separable and the convolution with it can be decomposed into two convolutions:

$C(s,t) * X(s,t) = CS(s) * CT(t) * X(s,t)$                   (8.3)

In that case the filters *CS(s)* and *CT(t)* can be inverted independently as a filter in one variable. Then we simply construct the inverse as a product of the two inverses. If we denote the inverses as *ZS(s)*- an inverse of *CS(s)*, and *ZT(t)*- the inverse of *CT(t)*, and the inverse of *C(s)* as *ZST(s,t)* we observe*:*

$ZST(s,t) = ZS(s) * ZT(t)$                          (8.4)

If *C(s, t)* is an exponent of a non-separable quadric with an *xy* term it can be brought to the main axis (diagonalized with a proper choice of a coordinate system of the *x, y* space) and in the system of coordinate *C(s,t)* would be separable.

## 9. Conclusions

We developed a method of decomposition of compactly supported symmetrical discreet filters into elementary filters. We show the criteria of invertibility for elementary filters. Thus, every filter can be presented as a composition of invertible and non-invertible component-filters. The effects of the invertible component can be fully undone. The non-invertible component leads to loss of resolution of the signal. In the space of a lesser resolution the signal can be presented without the effects of the instrument. The method of direct convolution inversion does not give any consideration for camera noise inevitable in applications and so we were pleasantly surprised that it can handle modest amount of camera noise. More studies are needed to explore the boundaries of applicability for diminishing signal to noise ratio.

## 10 Appendix – Examples and drawings

**Example 1. Some elementary (quadratic) filters can be further decomposed into "linear" filters.** Two notable examples of elementary non invertible convolutions.

When we presented a decomposition of symmetrical filters into elementary filters of 3 non-zero elements with quadratic characteristic polynomials, we were able to do it because every polynomial with real coefficients allows for a decomposition into polynomials of no more than second degree and quadratic polynomials are characteristic polynomials of elementary filters. However, if we are lucky, it may be possible to decompose a quadratic polynomial into a product of two linear polynomials with real coefficients.

Linear filters are nor symmetrical obviously (and neither are their linear char polynomials) but their convolutions may be. Let us consider two notable examples of such filters and polynomials.

**Example 1.1  *p=2***

A filter [*1, 2, 1*] (shown in Drawing 0. III) is in itself a convolutions of two equal filters of length 2:  [1, 1] and [1, 1]. Each of these two filters has a linear characteristic polynomial $p(x) = (x +1)$, and the product of the two characteristic polynomial equals the characteristic polynomial of the filter [1, 2, 1].

$$p(x) * p(x) = (x+1)^2 = x^2 + 2x + 1 \qquad (10.1)$$

In the case the characteristic polynomial has the root 1 with multiplicity 2.

Let's illustrate the duality between the linear filters and their kernels. Consider the filter [1, 1] – constant on final support of length 2 shown in Drawing 0. II and its kernel shown in Drawing 0. c) and d). Not surprisingly the kernel for both filters ([1, 1] and [1, 2, 1]) have the signal of the highest frequency wave (shown in Drawing 0 (c) and (d)).

**Example 1.2  *p= -2***

*Drawing 0.II* on the other hand shows a dual linear filter [-1, 1], with characteristic polynomial of *(x-1)* and its kernel – the constant signal. Note that [-1, 1] filter serve as an analogue of differentiation for the given discrete signals on the given lattice. It makes sense then that the "differentiation" maps constant signal to '0'. (It also maps linear signal to const, quadratic to liner etc.). Also note that the filter [-1, 1] convolved (composed) with itself is the other singular quadratic filter with $|p|=2$: [1, -2, 1] shown in Drawing 0.I

Its characteristic polynomial being

$$P(x) = (x-1) * (x-1) = x^2 - 2x = 1 \qquad (10.2)$$

***Drawing 0.***

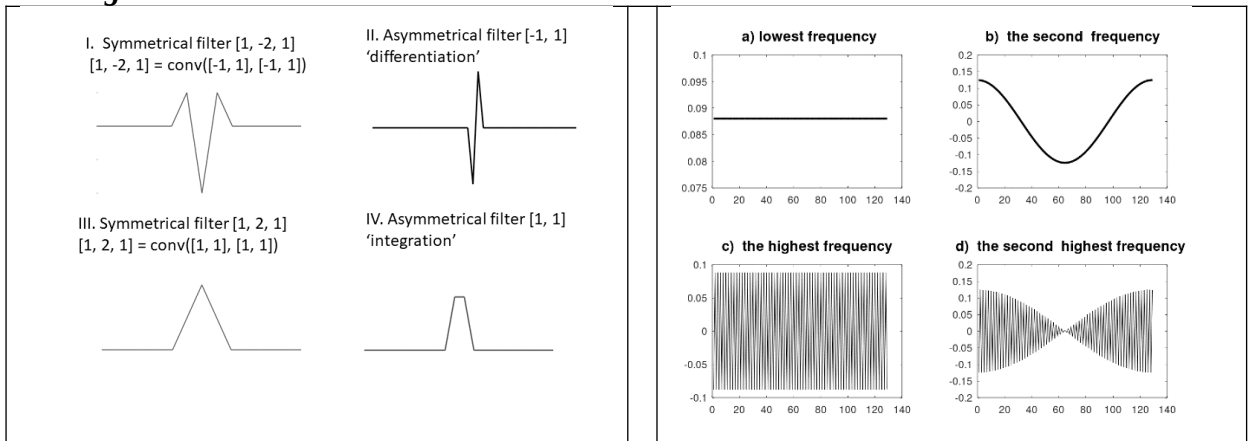

This filter is not invertible in any sense. Not in the sense of convergence to 0 (as t→∞), nor in a sense of limited values < Const (for any *t*) and being periodic (as t→∞).

We observe that the convolution with 'linear' filter *C0* (or with 'quadratic' filter *C0*C0*) maps the highest frequency of the given lattice *(0, 1, …,N)* to zero. That means that the signal Y - the convolution of the original signal *X* and the filter *C0*  *(Y =C0*X)*  does not have the highest frequency in its Fourier representation. Therefore, the Nyquist interval required to faithfully represent *Y* is larger than the Nyquist interval of the original signal *X*.

***Therefore, we are observing a loss of resolution in Y compared to resolution in X.***

Each of the linear (and non-symmetric) filters that we called 'differentiation' and 'integration' has one kernel function. The 'differentiation' filter *[1, -1]* zeros out constants and the' integration' *[1, 1]* zeros out the highest frequency. The corresponding quadratic filters of *[1, -2, 1]* (double differentiation) and *[1, 2, 1]* (double integration) have two kernel functions each. In addition to the highest frequency the integration zeros out the next highest. The kernel function though, looks pretty ugly, because the pixalation of *(1/N)* does not agree with the next highest frequency that would align with a pixelation of  *1/(N+1)* and so the kernel function looks ugly (see Drawing 0. d)). The double differentiation predictably maps constants and linear functions to zero.

**Example 2: Invertible filter** An elementary convolution E = [1, 2.3, 1] is invertible as it is an elementary symmetrical filter of the form [1, p, 1] with p =2.3 >2. Its inverse is shown in the Drawing 1a. Drawing 2a shows the distortion of the original signal caused by the filter E. The RMS error between the original and the distorted signal equals 0.47. The signal is restored with the direct inverse deconvolution, and the result is shown in Figure (b) of the Drawing 2. The RMS error is improved to 0.

**Example 3: Non-invertible filter**

Non-invertible filters allow for non-converging, but limited inverse filters. The inverses to [1, 1, 1] and [1, -1, 1] are shown respectively in the Figures (c) and (d) of the Drawing 1. The reconstruction from the effects of [1, 1, 1] and [1, -1, 1] are shown in Drawing 3a and (b)respectively. The reconstruction complete, meaning that the RMS error is zero, but the signal is shorter in length. The resolution of the signal (and its length), is reduced because of the existence of the kernel. There are two kernel functions associated with every non-invertible elementary symmetrical filter and therefore with each non-invertible elementary filter present in the decomposition of the original filter C the signal loses two units of length.

**Drawing 1**

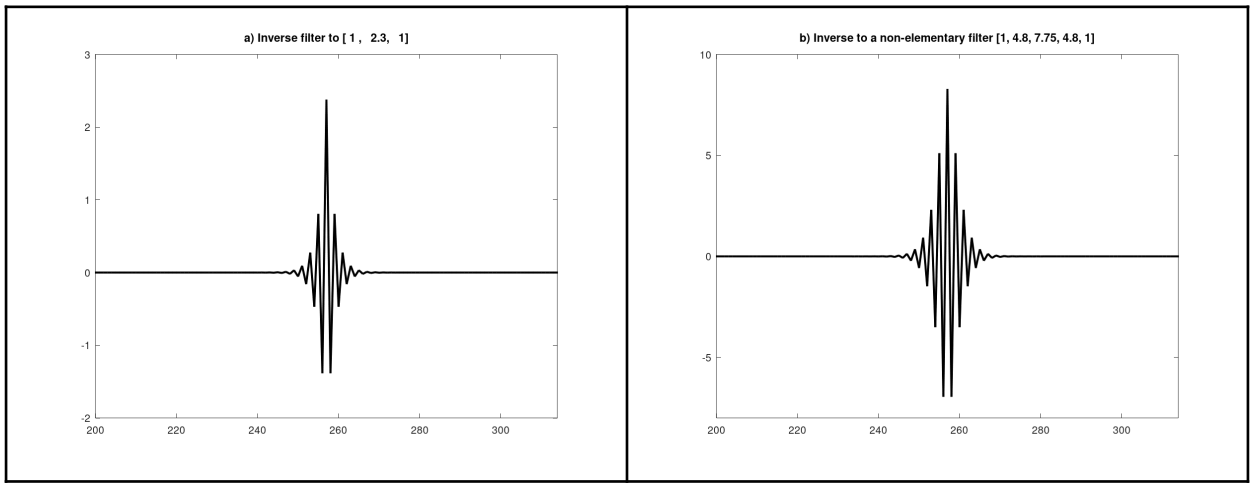

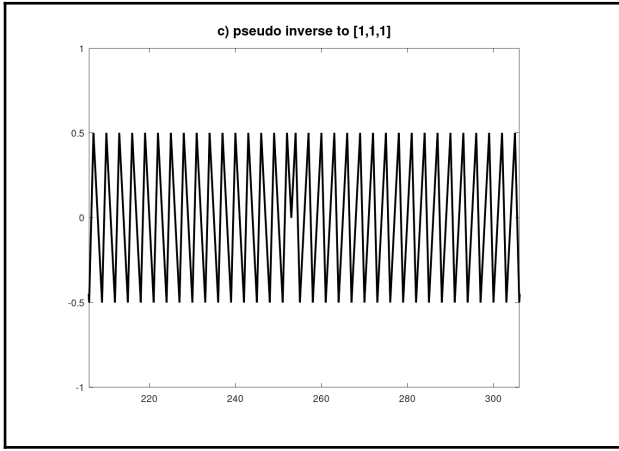
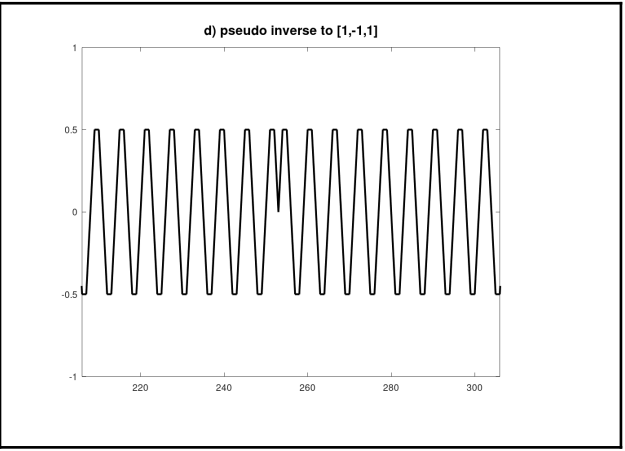

**Drawing 2**

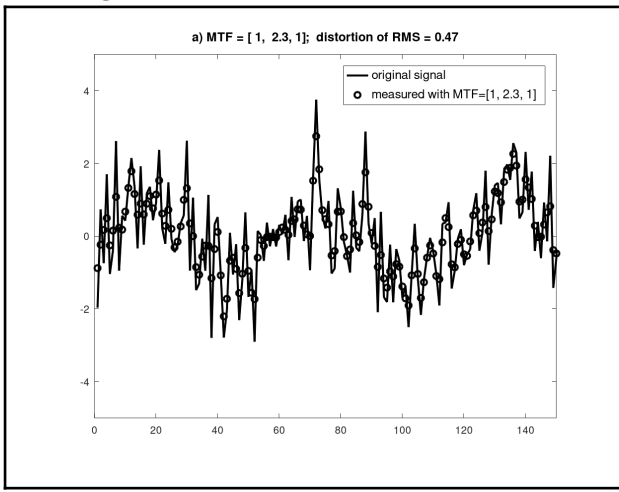
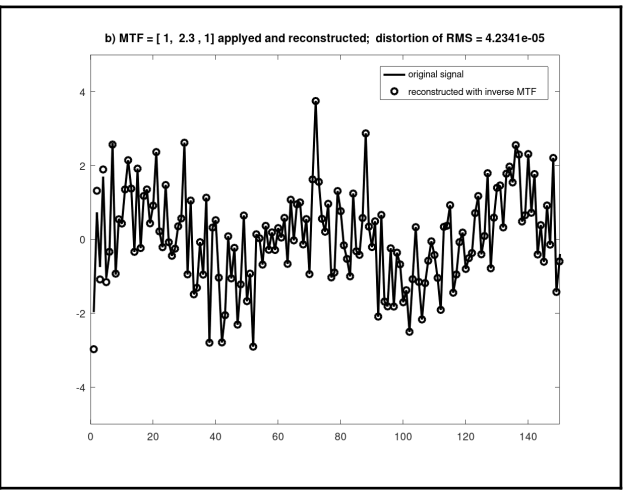

**Drawing 3**

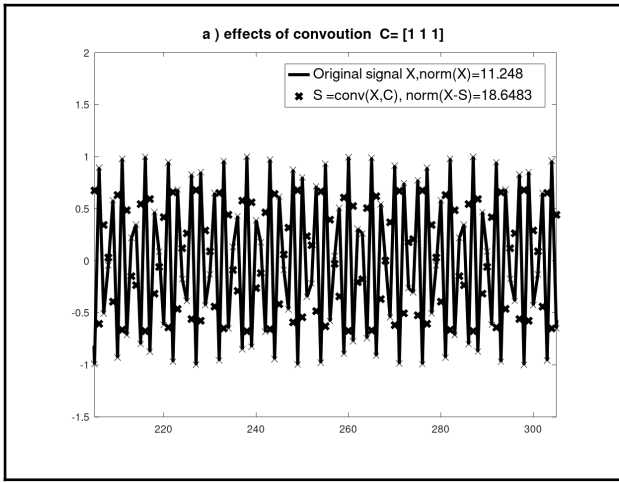
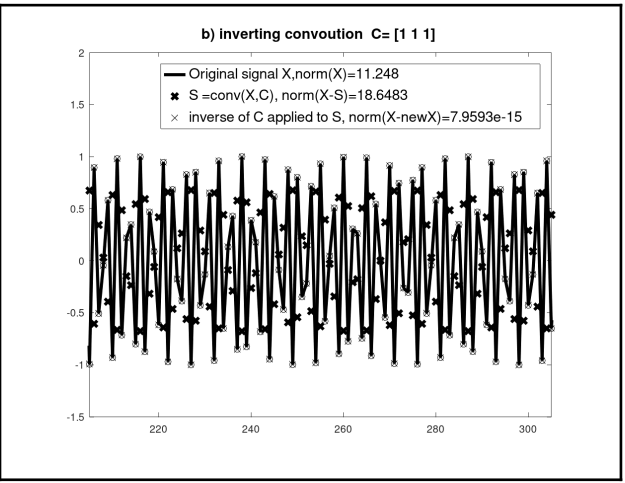

**Example 4**

Let us see how the direct inversion method compares with the industry-standard method for image deconvolution, namely, the Richardson Lucy (RL) deconvolution. RL Deconvolution is iterative and more computationally intensive than direct deconvolution disclosed in the present invention. Drawing 4 shows the checkerboard image (a), the checkerboard image **blurred with an invertible filter** (b), the image blurred reconstructed with direct inversion (c), and its RL reconstruction (d). Drawing 5 shows the blurred image with a significant amount of noise added to it. It appears that the noise overwhelms the RL method faster than it does the method of direct inversion. Drawing 6 (b)shows a blurred image A but no noise present. A Gaussian blur is used, and it is *a non-invertible filter*. Figure (c) shows a direct inverse deconvolution reconstruction of the image (a) and the reconstruction is almost complete and better than the RL reconstruction, shown in (d).

**Examples 5:** *Invertible 2D filter and its effects on Checkerboard image with camera noise added.* This example shows the effects of an invertible 2D filter and its inverse on the checkerboard image with some camera noise added. Drawings 5 shows the original checkerboard image (a), the image convolved with the invertible 2D filter (b), a 2-D reconstruction done with the direct inverse deconvolution (c), as compared with the performance of the industry-standard Richardson Lucy method shown in (d). Note that the method is not designed to handle noise, as it is a close term inversion of the operation of convolution. However, in practical applications, the 'camera' noise issue is ever present and it is important to assess whether the method offers any stability against it. The Drawings 5 and 7 illustrate that it does. The method can be combined with existing methods of denoising (such as median filtering) to increase its stability against the noise. Drawing 4 shows the checkerboard image (a), the checkerboard image **blurred with an invertible filter** (b), the image blurred reconstructed with direct inversion (c), and its RL reconstruction (d). Drawing 5 shows the blurred image with a significant amount of noise added to it.

**Drawing 5**                                    **Drawing 4**

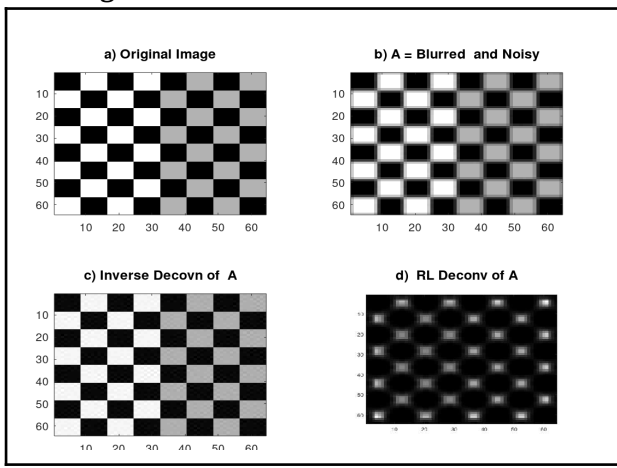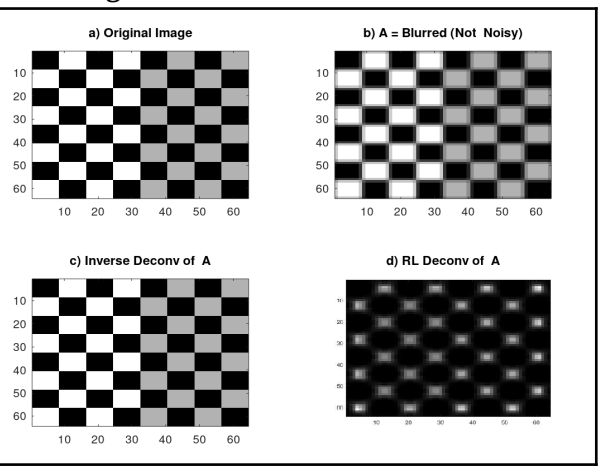

**Examples 6:** *Non-invertible 2D filter and its pseudo-inverse effects on Checkerboard image with no camera noise*

This example shows the effects of a non-invertible 2D filter and its pseudo-inverse on the checkerboard image without camera noise. As the blurring non-invertible Gaussian filter is used, the method uses a pseudo-inverse operation in combination with a slight reduction of the size of the image (image space factorization). Drawings 6 shows the original checkerboard image (a), the image convolved with the invertible 2D filter (b), a 2-D reconstruction done with the direct inverse deconvolution (c) as compared with the performance of the industry-standard Richardson Lucy method shown in (d). *No camera noise was added to the filtered image in this example.*

**Drawing 6**         **Drawing 7**

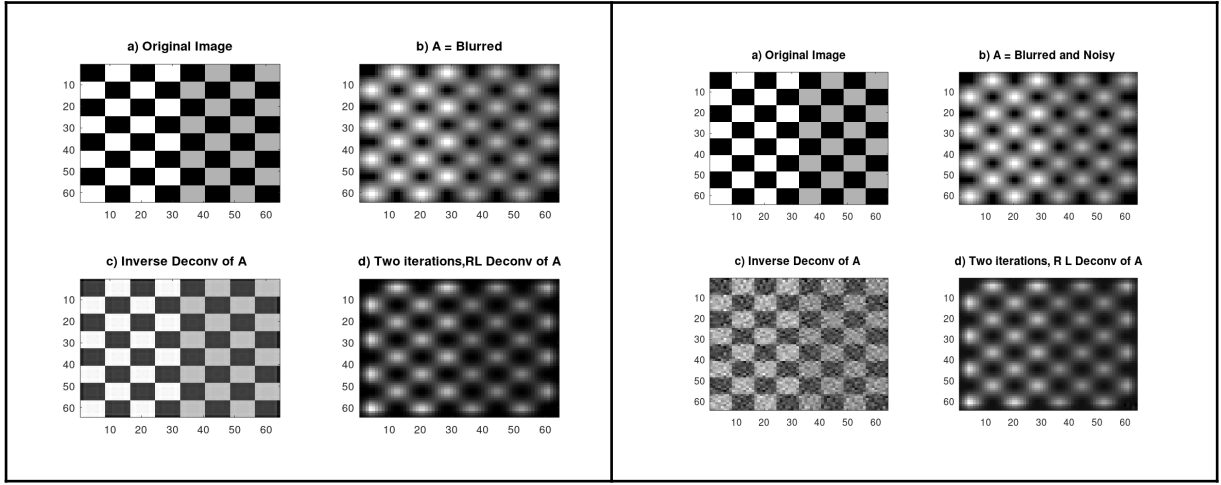

**Examples 7: Non-*invertible 2D filter and its pseudoinverse effects on Checkerboard image with camera noise added***

This example shows the effects of a non-invertible 2D filter and its pseudoinverse on the checkerboard image with some substantial amount of camera noise. As the blurring non-invertible Gaussian filter is used, the method uses a pseudo-inverse operation in combination with a slight reduction of the size of the image (image space factorization). Drawings 7 shows the original checkerboard image (a), the image convolved with the invertible 2D filter (b), a 2-D reconstruction done with the direct inverse deconvolution (c) as compared with the performance of the industry-standard Richardson Lucy method shown in (d). *Substantial amount of camera noise was added to the filtered image in this example.* One may observe that while a great deal of noise is present in the reconstruction done with the direct inversion, we clearly see the structure of the checkerboard, while it is damaged in the Richardson Lucy deconvolution. Arguably, the direct inverse deconvolution result is still useful because it reconstructs some of the features of the original image, such as the size of the checkerboards, better than the RL deconvolution does. The direct inversion method of deconvolution may perform better in a high noise situation if it is combined with the available methods of denoising (such as median filtering).